\documentclass{article}

\usepackage{amssymb,amsmath,bm,hyperref,geometry}
\usepackage{amssymb,amsmath,bm,graphics,graphicx,url,color}
\def\no{\noindent}
\def\pmatrix{\left(\begin{array}}
\def\endpmatrix{\end{array}\right)}
\def\dd{\mathrm{d}}

\newtheorem{theo}{Theorem}
\newtheorem{lem}{Lemma}
\newtheorem{cor}{Corollary}
\newtheorem{rem}{Remark}

\def\proof{\noindent\underline{Proof}\quad}
\def\QED{\mbox{~$\Box{~}$}}

\def\aa{\alpha}

\begin{document}

\title{A note on a stable algorithm for computing the fractional integrals of orthogonal polynomials}

\author{
P.\,Amodio\,\thanks{Dipartimento di Matematica, Universit\`a di Bari, Italy, {\tt pierluigi.amodio@uniba.it}} \and
L.\,Brugnano\,\thanks{Dipartimento di Matematica e Informatica ``U.\,Dini'', Universit\`a di Firenze, Italy, {\tt luigi.brugnano@unifi.it}} \and
F.\,Iavernaro\,\thanks{Dipartimento di Matematica, Universit\`a di Bari, Italy, {\tt felice.iavernaro@uniba.it}}}
%\\[.5cm] {\color{blue}\em\small --- Dedicated to John Butcher on the occasion of his 89th birthday ---}\\[.2cm]}

\maketitle

\begin{abstract} In this note we provide an algorithm for computing the fractional integrals of orthogonal polynomials, which is more stable than that using the expression of the polynomials w.r.t. the canonical basis. This algorithm is aimed at solving corresponding fractional differential equations. A few numerical examples are reported.

\bigskip
\no{\bf Keywords:} fractional differential equations, fractional integrals, orthogonal polynomials, Legendre polynomials, Chebyshev polynomials.

\bigskip
\no{\bf AMS-MSC:} 65L05, 65L03.

\end{abstract}

\section{Introduction} In recent years, the numerical solution of fractional differential equations has attracted the interest of many computational scientists, due to their usefulness in applications: we refer to the monograph \cite{Po1999} and to the review paper \cite{Ga2018} for some details. In this paper, we introduce some improvements to the implementation details related to the recent solution approach described in \cite{ABI2019} and based on previous work on HBVMs \cite{BIT2012,BS2014,BI2016,BI2018,ABI2020,ABI2022}, for solving fractional initial value problems in the form:
\begin{equation}\label{inifrac}
y^{(\aa)}(t) = f(t,y(t)), \qquad t\in[0,T], \qquad y(0)=y_0.
\end{equation}
Here, for ~$\aa\in(0,1]$, ~$y^{(\aa)}\equiv D^\aa y(t)$~ is the so-called Caputo fractional derivative:
\begin{equation}\label{Dalfa}
D^\aa g(t) = \frac{1}{\Gamma(1-\aa)}\int_0^t (t-x)^{-\aa} g'(x)\dd x.
\end{equation}
For sake of brevity, we assume that the equation is scalar and $\aa\in(0,1]$. However, the used argument can be easily generalized to values of $\aa>0$ and to systems of fractional differential equations, also having different orders. The Riemann-Liouville integral associated with (\ref{Dalfa}) is given by:
\begin{equation}\label{Ialfa}
I^\aa g(t) = \frac{1}{\Gamma(\aa)}\int_0^t (t-x)^{\aa-1} g(x)\dd x ~\equiv~ \frac{\Psi^\aa g(t)}{\Gamma(\aa)},
\end{equation}
where the function $\Psi^\aa g(t)$ has been defined for later use. 
It is known that \cite{Po1999}:
\begin{equation}\label{prop}
D^\aa I^\aa g(t) = g(t), ~\quad I^\aa D^\aa g(t)=g(t)-g(0), ~\quad
I^\aa t^j = \frac{j!}{\Gamma(\aa+j+1)}t^{j+\aa}, \quad j=0,1,2,\dots.
\end{equation}

In order to obtain a {\em quasi-polynomial} approximation of the solution of (\ref{inifrac}), according to \cite{ABI2020} we expand the vector field along an  orthonormal  polynomial basis,
\begin{equation}\label{orto}
\int_0^1 \omega(c)P_i(c)P_j(c)\dd c = \delta_{ij}, \qquad i,j=0,1,\dots,\qquad  
\end{equation}
with $\omega(c)$ a suitable weighting function,
\begin{equation}\label{wc}
\omega(c)\ge0, \qquad c\in[0,1], \qquad \int_0^1\omega(c)\dd c=1.
\end{equation}
One obtains
\begin{equation}\label{serie}
f(cT,y(cT)) = \sum_{j\ge0} \gamma_j(y) P_j(c), \qquad c\in[0,1],
\end{equation}
with the Fourier coefficients given by
\begin{equation}\label{gamj}
\gamma_j(y) = \int_0^1 \omega(c)P_j(c)f(cT,y(cT))\dd c, \qquad j=0,1,\dots.
\end{equation}
A polynomial approximation of degree $s-1$ to (\ref{serie}) is obtained by truncating the infinite series at the right-hand side after $s$ terms, which leads to the following fractional initial value problem approximating (\ref{inifrac}):
\begin{equation}\label{inifracs}
\sigma^{(\aa)}(cT) = \sum_{j=0}^{s-1} \gamma_j(\sigma) P_j(c), \qquad c\in[0,1], \qquad \sigma(0)=y_0,
\end{equation}
where $\gamma_j(\sigma)$ is formally given by (\ref{gamj}), upon replacing $y$ by $\sigma$. The solution of 
(\ref{inifracs}) is a {\em quasi-polynomial} of degree $s-1+\aa$, formally given by
\begin{equation}\label{sol}
\sigma(cT) = y_0 + h^\aa \sum_{j=0}^{s-1} \gamma_j(\sigma) I^\aa P_j(c), \qquad c\in[0,1].
\end{equation}
In order to compute the Fourier coefficients $\gamma_j(\sigma)$, one has to solve the system of equations:\footnote{For sake of brevity, hereafter we shall neglect the argument $\sigma$ in the Fourier coefficients.}
\begin{equation}\label{sist}
\gamma_j = \int_0^1 \omega(c)P_j(c)f\left( cT, y_0+h^\aa \sum_{k=0}^{s-1} \gamma_k I^\aa P_k(c)\right)\dd c, \qquad j=0,\dots,s-1.
\end{equation}
For this purpose, one may approximate the integrals in (\ref{sist}) with the corresponding Gaussian quadrature of order $2s$, whose abscissae are the zeros of $P_s$, $P_s(c_i)=0$, and with weights $b_i$, $i=1,\dots,s$, thus solving
\begin{equation}\label{sist1}
\gamma_j = \sum_{i=1}^s b_i P_j(c_i)f\left( c_iT, y_0+h^\aa \sum_{k=0}^{s-1} \gamma_k I^\aa P_k(c_i)\right), \qquad j=0,\dots,s-1.
\end{equation}
Consequently, one needs  to compute the fractional integrals:
\begin{equation}\label{Iaci}
I^\aa P_0(c_i), ~ I^\aa P_1(c_i), ~\dots,~ I^\aa P_{s-1}(c_i), \qquad i=1,\dots,s.
\end{equation}
The efficient and stable evaluation of the fractional integrals (\ref{Iaci}) is precisely the aim of this note.

\section{Computing the fractional integrals}
As is well-known, the family of orthonormal polynomials (\ref{orto})-(\ref{wc}) satisfies a 3-term recurrence
\begin{eqnarray}\nonumber
P_0(c) &\equiv& 1, \qquad\qquad\qquad\qquad\qquad\qquad\qquad c\in[0,1],\\
P_1(c) &=&(a_1c-b_1)P_0(c), \label{Pj}\\ \nonumber
P_j(c)  &=&(a_jc-b_j)P_{j-1}(c) - d_jP_{j-2}(c), \qquad j\ge 2,
\end{eqnarray} 
for suitable coefficients $a_j,b_j,d_j$, $j\ge1$ (with $d_1=0$).

As an example, the coefficients
\begin{equation}\label{ceby}
a_1=2, \quad b_1=1, \qquad a_j=4,\quad b_j=2,\quad d_j=1, \qquad j\ge2,
\end{equation}
provide the shifted and scaled Chebyshev polynomials of the first kind, corresponding to the choice of the weighting function ~ $\omega(c) = \left(\pi\sqrt{c(1-c)}\right)^{-1}$.

Another relevant example is provided by the shifted and scaled Legendre polynomials, corresponding to the weighting function ~$\omega(c)\equiv1$, for which:
\begin{equation}\label{legen}
b_j = \sqrt{4-j^{-2}} , \quad a_j = 2b_j, \quad j\ge 1,\qquad
d_j = \frac{j-1}j\sqrt{ \frac{2j+1}{2j-3} }, \quad j\ge 2.
\end{equation}
%\begin{eqnarray}\nonumber
%b_j &=& \sqrt{4-j^{-2}} , \qquad a_j = 2b_j, \qquad j\ge 1,\\ \label{legen}\\[-2mm]
%d_j &=& \frac{j-1}j\sqrt{ \frac{2j+1}{2j-3} }, \qquad j\ge 2.\nonumber
%\end{eqnarray}

 \smallskip
The following preliminary result holds true.

\begin{lem}\label{cgc} For any given $\aa>0$ and function $g(c)$, and with reference to the function $\Psi^\aa g(t)$ defined in (\ref{Ialfa}), one has: 
$$\Psi^\aa cg(c) = c\Psi^\aa g(c) -\Psi^{\aa+1}g(c).$$
\end{lem}
\proof In fact, from (\ref{Ialfa}) one has: 
\begin{eqnarray*}
\Psi^\aa cg(c) &=&\int_0^c (c-x)^{\aa-1}xg(x)\dd x
~=~\int_0^c (c-x)^{\aa-1}(x-c+c)g(x)\dd x\\[2mm]
&=& c\int_0^c (c-x)^{\aa-1}g(x)\dd x -\int_0^c (c-x)^\aa g(x)\dd x
~\equiv~ c \Psi^\aa g(c) -  \Psi^{\aa+1}g(c).\QED
\end{eqnarray*}
\smallskip

Consequently, from the last property in (\ref{prop}), Lemma~\ref{cgc}, and with reference to (\ref{Iaci}), the following results follows (we omit their straightforward proofs, for sake of brevity).

\begin{theo}\label{IaPj}
The fractional integrals of the polynomials (\ref{Pj}) are given by $$I^\aa P_j(c) = \Psi^\aa P_j(c)/ \Gamma(\aa),\qquad j\ge0,$$ where:
\begin{eqnarray}\nonumber
\Psi^\aa P_0(c) &=&\frac{c^\aa}{\aa},\\ \label{it3}
\Psi^\aa P_1(c) &=& (a_1c-b_1)\Psi^\aa P_0(c) - a_1 \Psi^{\aa+1}P_0(c),\\[2mm] \nonumber
\Psi^\aa P_j(c) &=& (a_jc-b_j)\Psi^\aa P_{j-1}(c) - a_j\Psi^{\aa+1}P_{j-1}(c) -d_j\Psi^\aa P_{j-2}(c), \qquad j\ge2.
\end{eqnarray}
\end{theo}

\smallskip
\begin{cor}\label{triango} For computing the fractional integrals
\begin{equation}\label{Iaj}
I^\aa P_0(c), ~ I^\aa P_1(c), ~\dots , ~I^\aa P_{s-2}(c), ~I^\aa P_{s-1}(c),
\end{equation}
by means of the recurrence scheme (\ref{it3}), one needs to compute the following triangular table:
\begin{equation}\label{table}
\begin{array}{llllll}
\Psi^\aa  P_0(c) \\
\Psi^{\aa+1} P_0(c) &\Psi^\aa  P_1(c)\\
\Psi^{\aa +2} P_0(c) &\Psi^{\aa +1} P_1(c) &\Psi^\aa  P_2(c)\\
\vdots                 &\vdots                 &\vdots          &\ddots \\ 
\Psi^{\aa +s-2} P_0(c) &\Psi^{\aa +s-3} P_1(c) &\Psi^{\aa +s-4} P_2(c) &\dots   &\Psi^\aa  P_{s-2}(c) \\
\Psi^{\aa +s-1} P_0(c)    &\Psi^{\aa +s-2} P_1(c) &\Psi^{\aa +s-3} P_2(c) &\dots   &\Psi^{\aa +1} P_{s-2}(c) &\Psi^\aa  P_{s-1}(c) \\
\end{array}\end{equation}
whose diagonal provides, upon division by $\Gamma(\aa)$, the integrals (\ref{Iaj}). 
\end{cor}

\begin{rem}\label{3vecs} Since (\ref{it3}) is a 3-term recurrence, in order to evaluate the integrals (\ref{Iaj}) only 3 vectors are actually needed, in place of the full triangular table (\ref{table}). For convenience, in  Appendix~\ref{code} we list a corresponding Matlab$^\copyright$ function.
\end{rem}

We conclude this section by recalling that, by using the known expressions of the Chebyshev and Legendre polynomials w.r.t. the canonical basis, and exploiting the last property in (\ref{prop}), one obtains
\begin{itemize}
\item for the Chebyshev polynomials:
\begin{eqnarray}\label{cebint}
I^\aa P_0(c) &=& \frac{c^\aa}{\Gamma(\aa+1)}, \\ \nonumber 
I^\aa P_j(c) &=& j\sqrt{2}\sum_{i=0}^j (-1)^{j-i}\frac{(j+i-1)! i!}{(j-i)!\Gamma(\aa+i+1) \prod_{k=1}^{2i}\frac{k}2 } c^{i+\aa}, \qquad j\ge1;
\end{eqnarray}

\item for the Legendre polynomials:
\begin{equation}\label{legint}
I^\aa P_j(c) = \sqrt{2j+1}\sum_{i=0}^j (-1)^{j-i}\frac{(j+i)!}{(j-i)! i! \Gamma(\aa+i+1)}c^{i+\aa}, \qquad j\ge0.
\end{equation}
Clearly, both (\ref{cebint}) and (\ref{legint}) can be efficiently evaluated by using a straightforward generalization of the Horner algorithm. In more detail, if
\begin{equation}\label{pol}
I^\aa P_j(c) = \sum_{i=0}^j p_i c^{i+\aa},
\end{equation}
then we can evaluate:\footnote{Of course, only one scalar $\rho$ is actually needed in (\ref{horn}): the subscripts have been added for sake of clarity.}
\begin{equation}\label{horn}
\rho_j = p_j, \qquad \rho_{i-1} = \rho_i c + p_{i-1}, \quad i= j,j-1,\dots,1, \qquad I^\aa P_j(c) \equiv \rho_0 c^\aa.
\end{equation}
\end{itemize}

\section{Numerical tests} 
In this section we report a few numerical tests illustrating the advantage of the proposed procedure (\ref{it3}) w.r.t. the standard one (\ref{horn}). All tests have been carried out on a 3GHz Intel Xeon W10 core computer with 64GB of memory, running Matlab$^\copyright$ 2020b.

At first, we show the maximum error in the computed fractional integrals (\ref{Iaj}), with $\aa=0.5$ and $s=25$, in the case of the Chebyshev and Legendre polynomials, respectively computed by using:
\begin{itemize}
\item the Horner algorithm (\ref{cebint})-(\ref{horn}) and the 3-term recurrence (\ref{ceby})-(\ref{it3});
\item  the Horner algorithm (\ref{legint})-(\ref{horn}) and the 3-term recurrence (\ref{legen})-(\ref{it3}). 
\end{itemize}
Reference integrals have been computed by using the variable precision arithmetic of the Matlab$^\copyright$ Symbolic Toolbox. The obtained results are plotted in Figure~\ref{fig3}: as one may see, in both cases, the error growth is much more favorable (and similar) when using the 3-term recurrence (\ref{it3}), w.r.t. the use of the canonical basis coupled with the Horner algorithm (\ref{horn}).

Next, we solve a fractional differential equation, using the expansions (\ref{inifracs})--(\ref{sist1}) of the vector field, along the Chebyshev basis (\ref{Pj})-(\ref{ceby}) and the Legendre basis (\ref{Pj})-(\ref{legen}). We compare the methods obtained by 
using the Horner algorithms (\ref{cebint})-(\ref{horn}) and   (\ref{legint})-(\ref{horn}), with the corresponding 3-term recurrences (\ref{ceby})-(\ref{it3}) and (\ref{legen})-(\ref{it3}), for computing the required fractional integrals (\ref{Iaci}).  We consider the following problem taken from \cite{Ga2018}:
\begin{eqnarray}\nonumber
y^{(\aa)}(t) &=& -y(t)^{\frac{3}2} +\frac{40320}{\Gamma(9-\aa)}t^{8-\aa} - 
3\frac{\Gamma(5+\frac{\aa}2)}{\Gamma(5-\frac{\aa}2)} t^{4-\frac{\aa}2} +
\left(\frac{3}2t^\frac{\aa}2-t^4\right)^3 +\frac{9}4\Gamma(\aa+1), \\ \label{prob1} \\[-3mm] 
&&t\in[0,T], \qquad y(0)=0,\nonumber
\end{eqnarray}
whose exact solution is known to be
$$y(t) = t^8-3t^{4+\frac{\aa}2}+\frac{9}4 t^\aa.$$
As emphasized in \cite{Ga2018}, this problem is surely of interest because, unlike several other problems often proposed in the literature, it does not present an artificial smooth solution, which is indeed not realistic in most of the fractional-order applications. We solve the problem for $\aa=0.5$, and $T=0.5$ and $T=1$. The results for the two cases are summarized in Figures~\ref{fig1} and \ref{fig2}, where we plot the error $\|y-\sigma\|$ w.r.t. the parameter $s$ in (\ref{inifracs}): 

\begin{figure}[hp]\centering
\includegraphics[width=7cm]{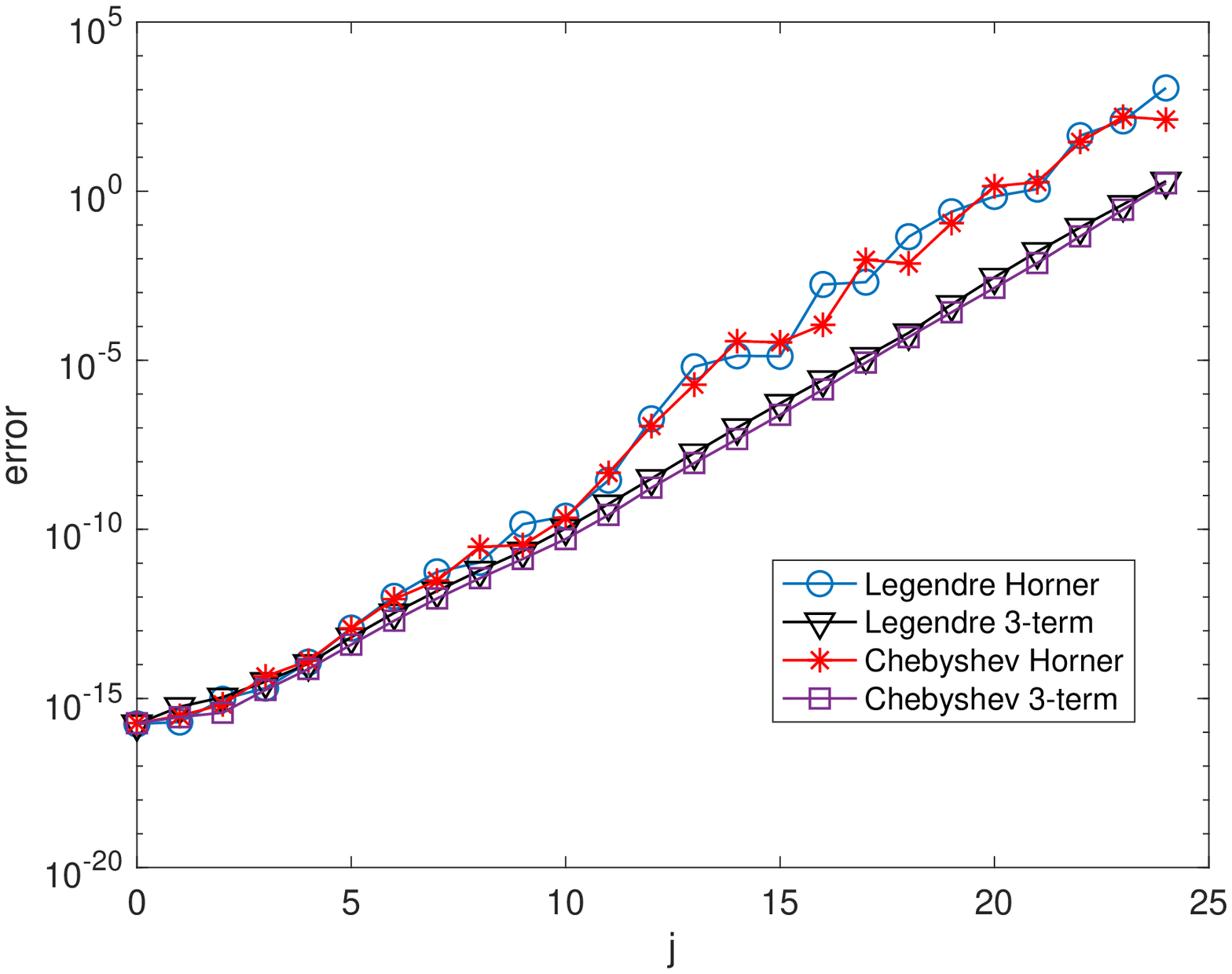}
\caption{Maximum error in the computed fractional integrals $I^\aa P_j(c)$, $\aa=0.5$, for the Chebyshev and Legendre polynomials, respectively computed by using the Horner algorithms (\ref{cebint})-(\ref{horn}) and (\ref{legint})-(\ref{horn}), and the 3-term recurrences  (\ref{ceby})-(\ref{it3}) and (\ref{legen})-(\ref{it3}).}
\label{fig3}

\bigskip
\bigskip
\includegraphics[width=4.8cm]{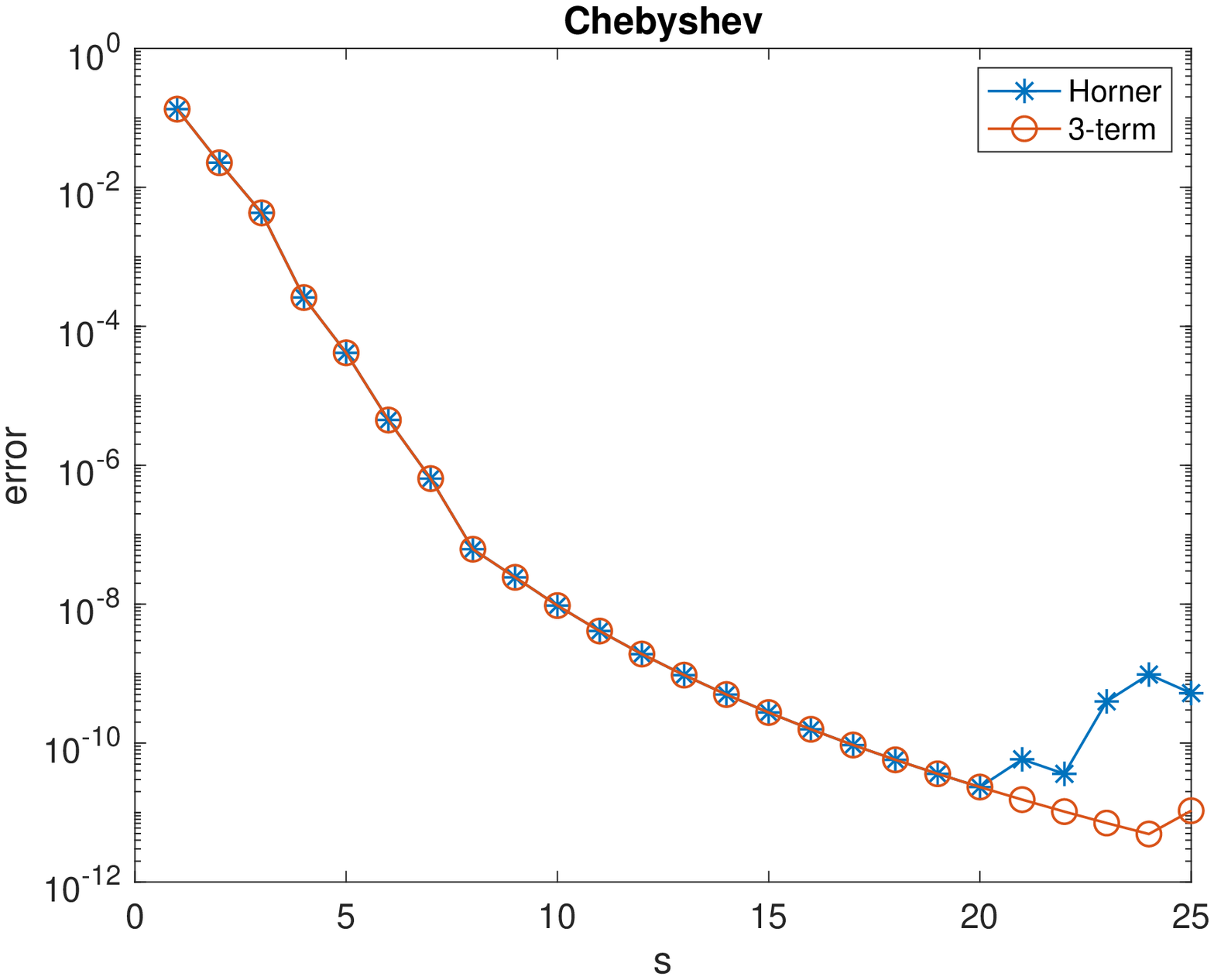}\quad
\includegraphics[width=4.8cm]{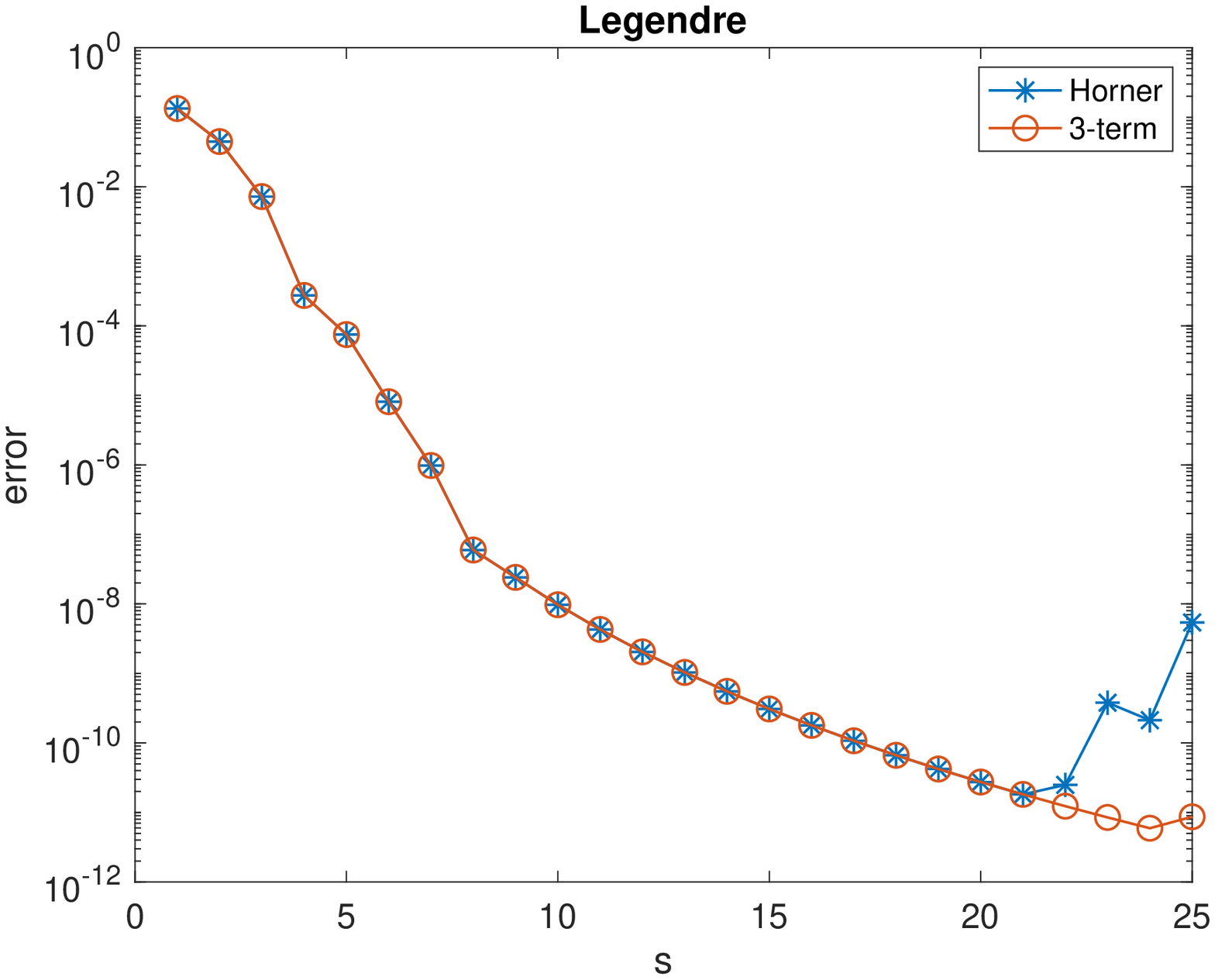}\quad
\includegraphics[width=4.8cm]{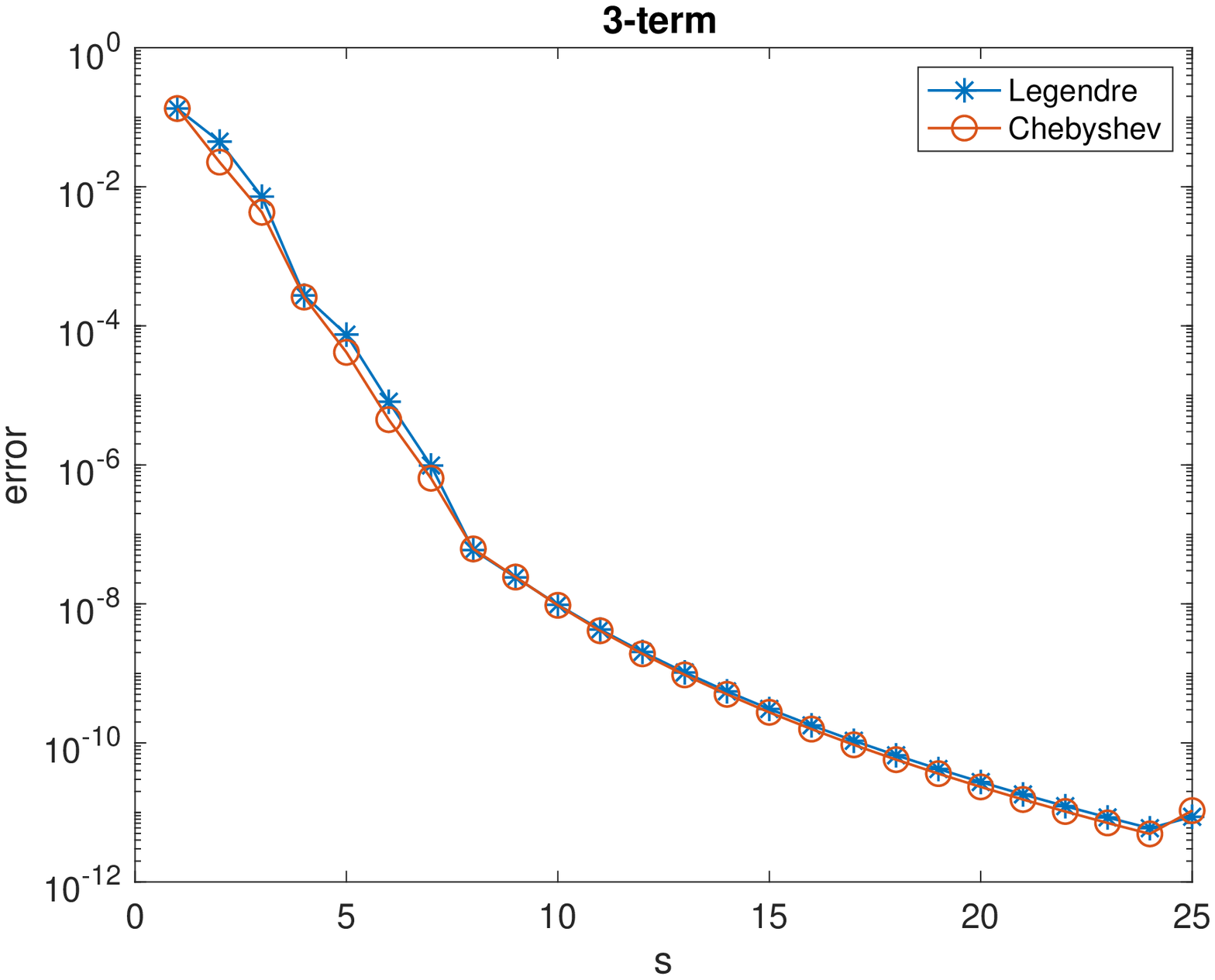}
\caption{Numerical solution of problem (\ref{prob1}), with $\alpha=0.5$ and $T=0.5$, by using the Chebyshev polynomials (left plot), the Legendre polynomials (middle plot), and both computed by using (\ref{it3}) (right plot).}
\label{fig1}

\bigskip
\bigskip
\includegraphics[width=4.8cm]{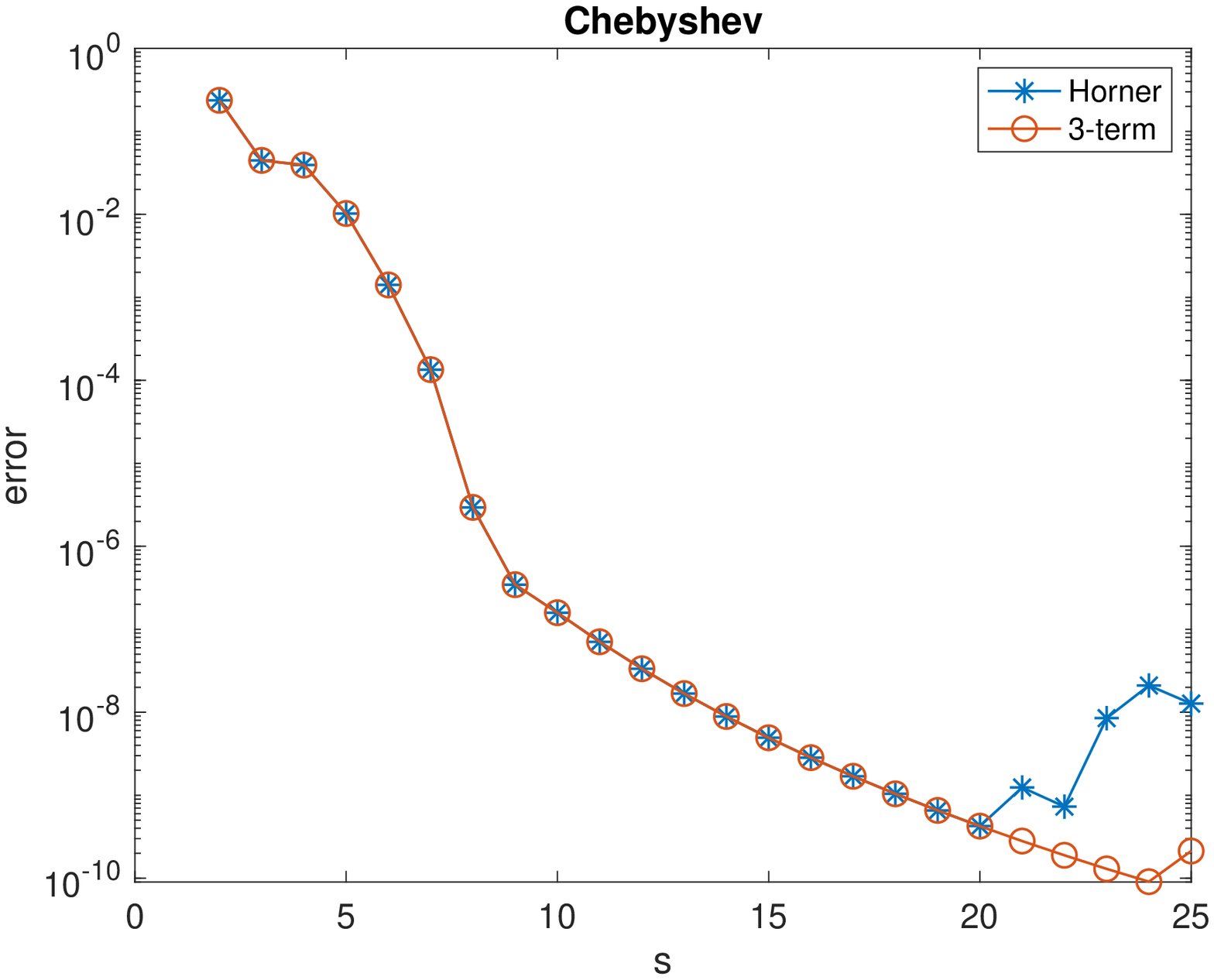}\quad
\includegraphics[width=4.8cm]{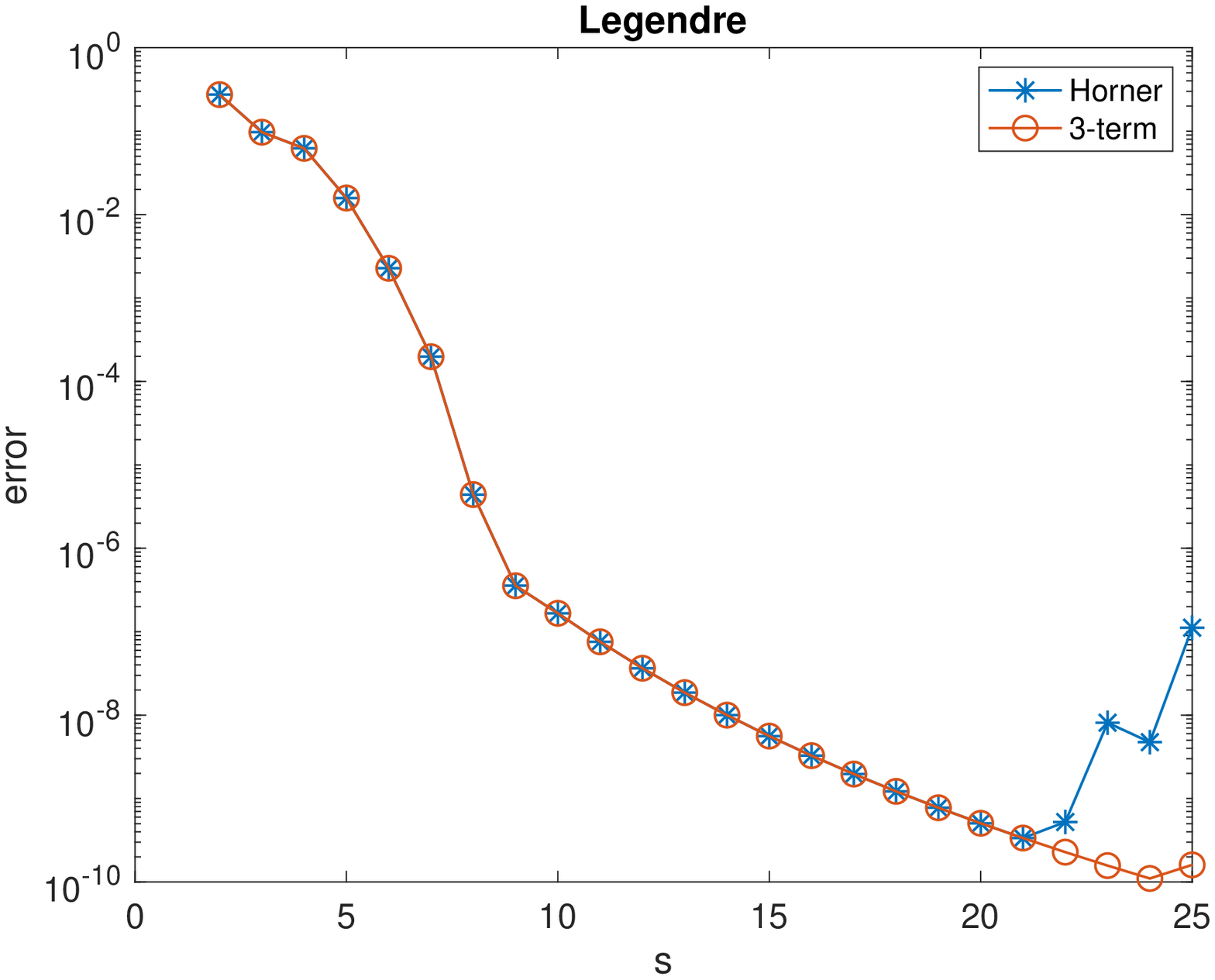}\quad
\includegraphics[width=4.8cm]{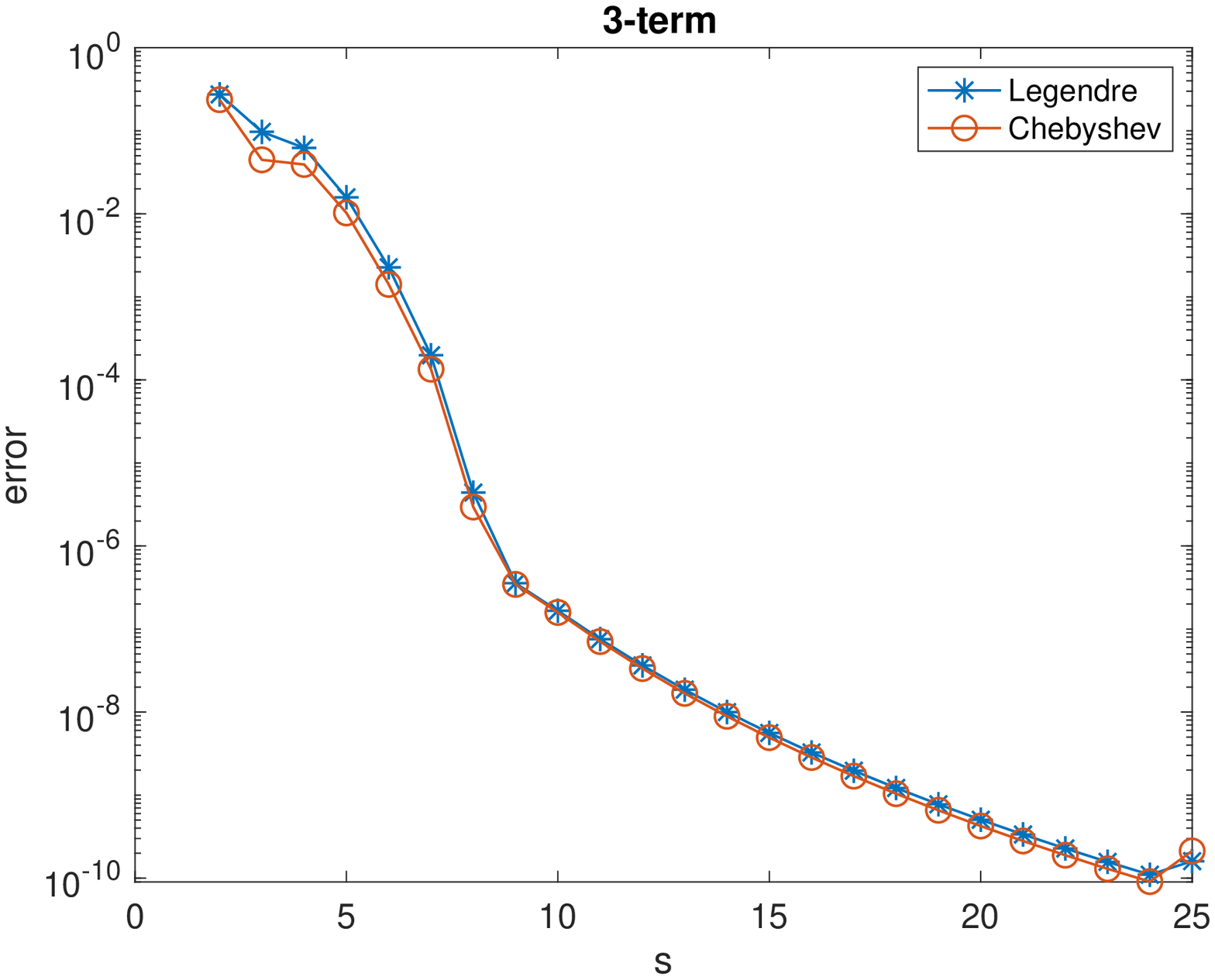}
\caption{Numerical solution of problem (\ref{prob1}), with $\alpha=0.5$ and $T=1$, by using the Chebyshev polynomials (left plot), the Legendre polynomials (middle plot), and both computed by using (\ref{it3}) (right plot).}
\label{fig2}
\end{figure}

\begin{itemize}
\item in both cases, as long as $s\le 20$, the errors produced by the 3-term recurrences (\ref{ceby})-(\ref{it3}) and (\ref{legen})-(\ref{it3})  is comparable (actually a bit smaller) than that obtained by the corresponding direct formulae (\ref{cebint}) and (\ref{legint}), computed via the Horner algorithm (\ref{horn}). A completely different behavior emerges for larger values of $s$. The error stagnates (actually, grows) when using the latter approach while, the better stability properties of the 3-term recurrence relations allows us to improve the accuracy of the approximation for values of $s$ up to $24$.

\item when using the 3-term recurrence (\ref{it3}), the approximation  errors produced by the Chebyshev and the Legendre polynomials are remarkably similar w.r.t. the parameter $s$.

\end{itemize}
\begin{rem}
It is worth noticing that the implementation of the fractional integrals (\ref{it3}) allows us to improve the results in \cite{ABI2019}, obtained by using the Legendre polynomials computed via (\ref{horn}). 
\end{rem}

\section{Conclusions} In this note, we have devised a 3-term recurrence for computing fractional integrals of orthogonal polynomials. The corresponding algorithm has been used for solving fractional differential equations, and proved to be more effective than the standard approach based on the canonical basis.

%\subsection*{Declarations} The authors declare no conflict of interests.

\paragraph*{Acknoledgements.} 
The authors wish to thank the {\em mrSIR} crowdfunding \cite{mrSIR} for the financial support.

\section{Appendix}\label{code}

\begin{verbatim}
function Ialfa = frac_int( a, b, d, alfa, c )
%
% Matlab function for computing the fractional integrals (18).
%
s        = length(a);
Ialfa    = zeros(s+1,1);
Psi1     = x.^(alfa+(0:s))./(alfa+(0:s)); 
Ialfa(1) = Psi1(1);
if s>=1
   Psi2     = ( a(1)*x-b(1) )*Psi1(1:s) - a(1)*Psi1(2:s+1);
   Ialfa(2) = Psi2(1); 
   for j = 2:s   
       Psi0       = Psi1; 
       Psi1       = Psi2;
       Psi2       = ( a(j)*x-b(j) )*Psi1(1:s-j+1) - a(j)*Psi1(2:s-j+2) ...
                                                  - c(j)*Psi0(1:s-j+1);
       Ialfa(j+1) = Psi2(1); 
   end
end   
Ialfa = Ialfa/gamma(alfa);
return
\end{verbatim}


\begin{thebibliography}{9}

\bibitem{ABI2019} P.\,Amodio, L.\,Brugnano, F.\,Iavernaro. Spectrally accurate solutions of nonlinear fractional initial value problems. {\em AIP Conf. Proc.} {\bf 2116} (2019) 140005. \url{https://doi.org/10.1063/1.5114132}

\bibitem{ABI2020} P.\,Amodio, L.\,Brugnano, F.\,Iavernaro. Analysis of Spectral Hamiltonian Boundary Value Methods (SHBVMs) for the numerical solution of ODE problems.  {\em Numer. Algorithms}   {\bf 83} (2020) 1489--1508. \url{https://doi.org/10.1007/s11075-019-00733-7}

\bibitem{ABI2022} P.\,Amodio, L.\,Brugnano, F.\,Iavernaro. Arbitrarily high-order energy-conserving methods for Poisson problems. {\em Numer. Algoritms}  (2022). \url{https://doi.org/10.1007/s11075-022-01285-z}

\bibitem{BI2016} L.\,Brugnano, F.\,Iavernaro. {\em Line Integral Methods for Conservative Problems}.  Chapman et Hall/CRC, Boca~Raton, FL, USA, 2016.

\bibitem{BI2018}  L.\,Brugnano, F.\,Iavernaro. Line Integral Solution of Differential Problems. {\em Axioms} {\bf 7}(2) (2018) 36. \url{https://doi.org/10.3390/axioms7020036}

\bibitem{BIT2012} L.\,Brugnano, F.\,Iavernaro, D.\,Trigiante.  A simple framework for the derivation and analysis of effective one-step methods for ODEs. {\em Appl. Math. Comput.} {\bf 218} (2012) 8475--8485. \url{https://doi.org/10.1016/j.amc.2012.01.074}

\bibitem{BS2014} L.\,Brugnano, Y.\,Sun. Multiple invariants conserving Runge-Kutta type methods for Hamiltonian problems. {\em Numer. Algorithms} {\bf 65} (2014) 611--632. \url{https://doi.org/10.1007/s11075-013-9769-9}

%\bibitem{CoNa2001} F.\,Costabile, A.\,Napoli. A method for global approximation of the initial value problem. {\em Numer. Algorithms} {\bf 27 } (2001) 119--130.

\bibitem{Ga2018} R.\,Garrappa. Numerical solution of fractional differential equations: a survey and a software tutorial. {\em Mathematics} {\bf 6}(2) (2018) 16. \url{http://doi.org/10.3390/math6020016}

\bibitem{Po1999} I.\,Podlubny. {\em Fractional differential equations. An introduction to fractional derivatives, fractional differential equations, to methods of their solution and some of their applications.}  Academic Press, Inc., San Diego, CA, 1999.

\bibitem{mrSIR} \url{https://www.mrsir.it/en/about-us/}


\end{thebibliography}
\end{document}